\documentclass[10pt]{amsart}
\usepackage{amssymb}
\usepackage{enumerate}
\usepackage{epsf}
\usepackage{psfrag}
\DeclareGraphicsExtensions{.eps,.art,.ART,.ps}
\newcommand{\bbR}{\mathbb{R}}      

\begin{document}
\begin{abstract}
We give an elementary derivation of the Montgomery phase formula for the motion of an Euler top, using only basic facts about the Euler equation and parallel transport on the $2$-sphere (whose holonomy is seen to be responsible for the geometric phase). We also give an approximate geometric interpretation of the geometric phase for motions starting close to an unstable equilibrium point.
\end{abstract}
%
%
\title[An elementary derivation of the Montgomery phase formula]{An elementary derivation of the Montgomery phase formula for the Euler top}
\author{Jos\'{e} Nat\'{a}rio}
\address{Centro de An\'alise Matem\'atica, Geometria e Sistemas Din\^amicos, Departamento de Matem\'atica, Instituto Superior T\'ecnico, 1049-001 Lisboa, Portugal}
\thanks{Partially supported by FCT (Portugal).}
\maketitle
%
%
%
\section*{Introduction}
The motion of the Euler top is governed by the Euler equation, whose generic orbits are periodic. By applying the Stokes theorem to a suitable surface on $T^*SO(3)$, Montgomery \cite{Montgomery91} obtained a Berry-Hannay-like formula for the angle by which the final position of the Euler top is rotated with respect to the initial position after one period (see also \cite{MMR90, Marsden92}).  A different derivation, based on the Poinsot description of the motion and the Gauss-Bonnet theorem, was given by Levi \cite{Levi93}.\footnote{Both these derivations have been generalized to more complicated mechanical systems (see \cite{AKS95, Cabrera07}).}

The purpose of the present paper is to give a third, more elementary derivation, in the sense that it utilizes only basic facts about the Euler equation and parallel transport on the $2$-sphere. The basic observation is that the motion of a fixed orthonormal basis as seen in the Euler top's frame can be easily understood in terms of the Euler flow on a sphere of fixed angular momentum norm and parallel transport on this sphere.

The structure of the paper is as follows: the first section briefly reviews the theory of the Euler top; the main result is proved in the second section, with Montgomery's formula deduced as a corollary in the third section; the fourth section contains an approximate geometric interpretation of the geometric phase for motions starting close to an unstable equilibrium point.
%
%
%
\section{Euler top}
In this section we briefly review the theory of the Euler top, mainly to fix the notation. This material is standard and can be found in almost any book on mechanics (e.g.~\cite{Ar97, Goldstein02, MR99, O02}).

 An {\bf Euler top} is a rigid body with a fixed point moving freely in an inertial frame. Its motion is described by a curve $S:\bbR \to SO(3)$ which at each instant gives the orientation of the body with respect to a reference position. We have $\dot{S}(t)=S(t)A(t)$, where $A(t) \in \mathfrak{so}(3)$ can be identified with ${\bf \Omega}(t) \in \bbR^3$ through the well known linear isomorphism defined by
\[
A(t) {\bf \xi} = {\bf \Omega}(t) \times {\bf \xi}
\]
for all ${\bf \xi} \in \bbR^3$. It is then easily shown that $\boldsymbol{\omega}(t)=S(t){\bf \Omega}(t)$ is the angular velocity of the Euler top, and hence ${\bf \Omega}(t)$ is the angular velocity as seen in the Euler top's frame.

The (conserved) total angular momentum of the Euler top is given by
\[
{\bf p} = S(t)I{\bf \Omega}(t),
\]
where $I:\bbR^3 \to \bbR^3$ is a symmetric, positive definite linear operator whose eigenvalues $I_1,I_2,I_3$ (called the {\bf principal moments of inertia}) satisfy the triangular inequality $I_1 \leq I_2 + I_3$ (and cyclic permutations). The $1$-dimensional eigenspaces of $I$ are called the {\bf principal axes of inertia}.

\begin{figure}[h!]
\begin{center}
\psfrag{e1}{$e_1$}
\psfrag{e2}{$e_2$}
\psfrag{e3}{$e_3$}
\epsfxsize=.7\textwidth
\leavevmode
\epsfbox{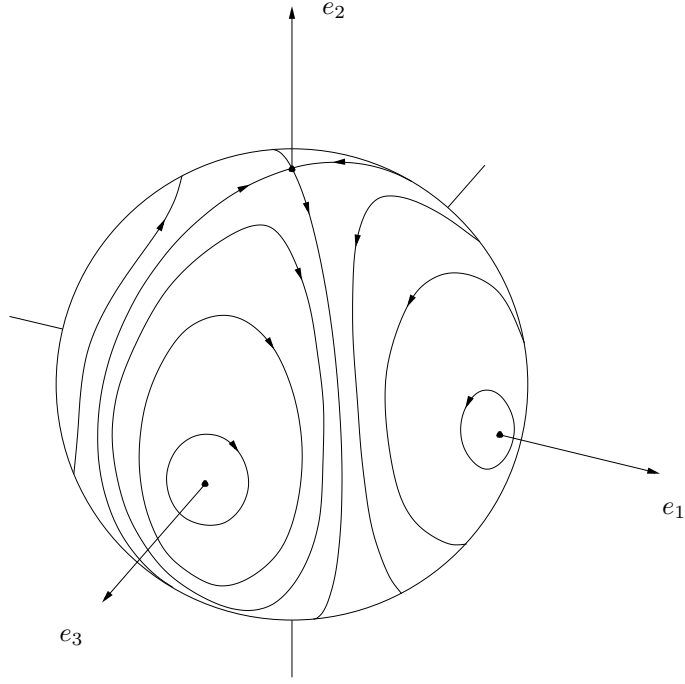}
\end{center}
\caption{Integral curves of the Euler equations on the angular momentum sphere, assuming $I_1>I_2>I_3$.}\label{Figure1}
\end{figure}

The vector ${\bf P}(t)=S^{-1}(t){\bf p}$ is the angular momentum vector as seen in the Euler top's frame, and is in general not constant. In fact,
\[
{\bf{\dot p}} = {\bf 0} \Leftrightarrow \dot{S}{\bf P}+S{\bf{\dot P}}={\bf 0}\Leftrightarrow S({\bf \Omega}\times{\bf P}) + S{\bf{\dot P}}={\bf 0},
\]
and hence ${\bf P}(t)$ satisfies the so-called {\bf Euler equation}
\[
{\bf{\dot P}}+{\bf \Omega}\times{\bf P}={\bf 0}.
\]
Using ${\bf \Omega}=I^{-1}{\bf P}$, this is a ODE for ${\bf P}(t)$, which in an orthonormal basis $\{{\bf e}_1, {\bf e}_2, {\bf e}_3 \}$ of eigenvectors of $I$ is written
\[
\begin{cases}
\displaystyle \dot{P}_1 = \left( \frac1{I_3} - \frac1{I_2}\right) P_2 P_3 \\
\displaystyle \dot{P}_2 = \left( \frac1{I_1} - \frac1{I_3}\right) P_3 P_1 \\
\displaystyle \dot{P}_3 = \left( \frac1{I_2} - \frac1{I_1}\right) P_1 P_2 
\end{cases}
\]
The Euler equations are integrable. An obvious first integral is $\|{\bf P}\|=\|{\bf p}\|=p$ (which we assume to be nonzero), and a second first integral is the kinetic energy $K=\frac12\langle{\bf P},{\bf \Omega}\rangle$:
\[
\frac{dK}{dt}=\frac12\frac{d}{dt}\langle{\bf P},I^{-1}{\bf P}\rangle = \langle{\bf{\dot P}},I^{-1}{\bf P}\rangle = \langle{\bf{\dot P}},{\bf \Omega}\rangle = - \langle{\bf \Omega}\times{\bf P},{\bf \Omega}\rangle = 0.
\]
Thus the solutions with fixed norm of the angular momentum live on the intersections of the angular momentum sphere $\|{\bf P}\|=p$ with the ellipsoids $\langle{\bf P},I^{-1}{\bf P}\rangle=2K$ for different values of $K$ (cf.~Figure~\ref{Figure1}).\footnote{In fact, it is easily shown that the Euler equations give the Hamiltonian flow of $K$ on the sphere $\|{\bf P}\| = p$ for the symplectic structure given by $\frac1p$ times the standard surface area element.}

%
%
%
\section{A picture of the motion}
Consider a right-handed orthonormal basis $\{ {\bf e}, {\bf f}, {\bf n} \}$, fixed in the inertial frame, with ${\bf n} = {\bf p}/p$. Then ${\bf N}(t)=S^{-1}(t)\,{\bf n}={\bf P}(t)/p$ will be the unit normal vector to the angular momentum sphere at the point ${\bf P}(t)$, and ${\bf E}(t)=S^{-1}(t)\,{\bf e}$ and ${\bf F}(t)=S^{-1}(t)\,{\bf f}$ will be (orthogonal) unit vectors tangent to the sphere at the same point.

We can ask ourselves what is the covariant derivative of ${\bf E}(t)$ as it moves along the sphere. Since ${\bf e}$ is constant, we have the Euler-like equation
\[
{\bf \dot{E}}+{\bf \Omega}\times{\bf E}={\bf 0}.
\]
As is well known \cite{Boothby03, Carmo76, Carmo93, ONeill83}, the covariant derivative is just the projection ${\bf \dot{E}}^\top$ of ${\bf \dot{E}}$ on the tangent space to the sphere. Thus\footnote{Here we use the vector identities $\langle {\bf a}, {\bf b} \times {\bf c} \rangle = \langle {\bf b}, {\bf c} \times {\bf a} \rangle$ and ${\bf a} \times ({\bf b} \times {\bf c}) = \langle {\bf a}, {\bf c} \rangle \, {\bf b} - \langle {\bf a}, {\bf b} \rangle \, {\bf c}$.}
\begin{align*}
{\bf \dot{E}}^\top & = \langle {\bf \dot{E}}, {\bf E} \rangle \, {\bf E} + \langle {\bf \dot{E}}, {\bf F} \rangle \, {\bf F} = - \langle {\bf \Omega}\times{\bf E}, {\bf E} \rangle \, {\bf E} - \langle {\bf \Omega}\times{\bf E}, {\bf F} \rangle \, {\bf F} \\
& = - \langle {\bf \Omega} \times{\bf E}, {\bf N} \times {\bf E} \rangle \, {\bf F} =  - \langle {\bf N}, {\bf E} \times ({\bf \Omega} \times{\bf E}) \rangle \, {\bf F} \\
& = - \langle {\bf N}, \|{\bf E}\|^2 {\bf \Omega} - \langle{\bf E},{\bf \Omega}\rangle \, {\bf E} \rangle \, {\bf F} = - \langle {\bf N}, {\bf \Omega} \rangle \, {\bf F} \\
& = - \frac1p\langle {\bf P}, {\bf \Omega} \rangle \, {\bf F} = - \frac{2K}p\,{\bf F}.
\end{align*}
In other words, ${\bf E}(t)$ rotates with constant angular velocity $-2K/p$ with respect to a parallel-transported frame on the sphere (and the same must obviously be true for ${\bf F}(t)$). This gives a nice picture of the motion: the fixed orthonormal basis $\{ {\bf e}, {\bf f}, {\bf n} \}$ is seen in the Euler top's frame as the moving orthonormal basis $\{ {\bf E}(t), {\bf F}(t), {\bf N}(t) \}$. If we imagine this basis placed on the point ${\bf P}(t)$ of the angular momentum sphere (which moves according to the Euler equation), then ${\bf N}(t)$ is the outward unit normal vector and ${\bf E}(t)$, ${\bf F}(t)$ rotate at a constant rate about ${\bf N}(t)$.
%
%
%
\section{Montgomery's formula}
Except for the six equilibrium points and four heteroclinic orbits, all orbits of ${\bf P}(t)$ on the angular momentum sphere are periodic (cf.~Figure~\ref{Figure1}). After one period, the angular momentum vector as seen in the Euler top's frame has returned to the initial position. Therefore the Euler top has rotated by some angle $\Delta \theta$ around {\bf p} with respect to its initial orientation. Montgomery's formula gives this angle as 
\[
\Delta \theta = \frac{2KT}{p} - \frac{A}{p^2},
\]
where $T$ is the period and $A$ is the area of the region bounded by the orbit on the angular momentum sphere (with sign according to whether the orientation of the periodic orbit is induced by the standard orientation of this region).\footnote{The two possible choices, $A$ and $-(4\pi p^2-A)=A - 4\pi p^2$, lead to the same angle $\text{mod}\, 2 \pi$. It is possible to make sense of the {\bf total} angle as $2\pi$ times the rotation number on the corresponding invariant torus, but its computation requires solving the Euler equations \cite{BCS05}.} The first term, which is just the time integral of the component of $\boldsymbol{\omega}(t)$ along ${\bf p}$, is a dynamical phase. The second term, which can be thought of as a solid angle, is what is usually called a geometric phase.

We can use our picture of the motion to give an elementary derivation of the Montgomery formula. First, it is clear that $\Delta \theta$ exactly coincides with minus the angle by which ${\bf E}(t)$ has rotated with respect to its initial position after one orbit. Since ${\bf E}(t)$ rotates with constant angular velocity $-2K/p$ with respect to a parallel-transported frame on the sphere, it will have rotated $-2KT/p$ with respect to such a frame after one orbit (dynamical phase). By elementary geometry of curved surfaces \cite{Carmo76, Carmo94, Morgan98, Spivak3}, upon returning to the initial point the parallel-transported frame will be rotated by an angle equal to the curvature of the sphere ($1/p^2$) times the (signed) area $A$ of the region enclosed by its path (geometric phase). Therefore we obtain
\[
\Delta \theta = - \left( - \frac{2KT}{p} + \frac{A}{p^2} \right) = \frac{2KT}{p} - \frac{A}{p^2},
\]
which is Montgomery's formula. In particular, the geometric phase is seen as a consequence of parallel transport on the angular momentum sphere.
%
%
%
\section{Interpretation of the geometric phase}
It is possible to give an approximate geometric interpretation of the geometric phase in the case of an orbit which starts close to one of the unstable equilibrium points on the angular momentum sphere (cf.~Figure~\ref{Figure1}). For such orbits, the motion of the Euler top is roughly as follows: initially it is rotating with angular speed $2K/p \simeq p/I_2$ about ${\bf p}$ (that is, ${\bf P}(t)$ remains close to the equilibrium point); then it quickly reorients itself by rotating $180^{\circ}$ about an axis perpendicular to ${\bf p}$ (${\bf P}(t)$ rapidly moves to the other unstable equilibrium point along a heteroclinic);\footnote{This is the twisting phenomenon studied in \cite{AAC91, CB97}.} in this upside-down orientation it continues to rotate about {\bf p} with the same angular speed $2K/p$ (${\bf P}(t)$ remains close to the other equilibrium point); finally, it performs another reorientation by rotating $180^{\circ}$ about a (different) axis perpendicular to ${\bf p}$ (${\bf P}(t)$ returns to the first equilibrium point along another heteroclinic). In other words, with respect to an auxiliary frame in uniform rotation about ${\bf p}$ with angular speed $2K/p \simeq p/I_2$, the motion of the Euler top during one period consists basically of two $180^\circ$ rotations about certain axes perpendicular to ${\bf p}$.\footnote{More precisely, if $\{ {\bf g}, {\bf h}, {\bf n} \}$ is an orthonormal basis fixed in the auxiliary frame, so that ${\bf{\dot g}} = \frac{2K}{p^2} {\bf p}\times {\bf g}$, then ${\bf G}=S^{-1}{\bf g}$ satisfies ${\bf{\dot G}}^\top = {\bf 0}$ (and similarly ${\bf{\dot H}}^\top = {\bf 0}$).} Since the auxiliary frame rotates by $2KT/p$ during one period, the geometric phase must result from the composition of the two reorientations (which in particular has to be a rotation about ${\bf p}$).

It is easily seen that the heteroclinics are half great circles.\footnote{Assuming $I_1>I_2>I_3$ they are contained in the planes $P_3 = \pm \sqrt{\frac{(I_1-I_2) I_3}{I_1 (I_2 - I_3)}}\,P_1$.} During the reorientations, the auxiliary frame rotates with respect to the Euler top's frame about the axes perpendicular to these great circles, in the direction given by the Euler flow. According to the geometric construction described in the appendix, the composition of these two rotations is a rotation about the middle moment of inertia axis by twice the angle between the two heteroclinics. But this is exactly $A/p^2$, where $A$ is the (signed) area enclosed by the two heteroclinics. Therefore the Euler top rotates by $-A/p^2$ about ${\bf p}$ with respect to the auxiliary frame, giving the geometric phase. 
%
%
%
\section*{Appendix}
Here we describe the geometric construction for the composition of rotations given in \cite{Penrose04}. This construction is dual to the more familiar Rodrigues-Hamilton construction, which uses triangles with vertices on the rotation axes \cite{Needham97, Pars65, Stillwell92}.

\begin{figure}[h!]
\begin{center}
\psfrag{t}{$\frac{\theta}2$}
\psfrag{f}{$\frac{\varphi}2$}
\psfrag{p}{$\frac{\psi}2$}
\psfrag{O}{$O$}
\psfrag{A}{$A$}
\psfrag{B}{$B$}
\epsfxsize=.5\textwidth
\leavevmode
\epsfbox{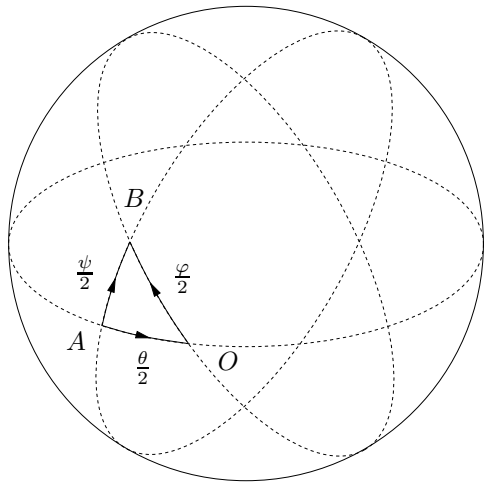}
\end{center}
\caption{Geometric construction for the composition of rotations.}\label{Figure2}
\end{figure}

Suppose that one wants to compose a rotation about some axis $a$ by an angle $\theta$ with a rotation about some other axis $b$ by an angle $\varphi$. Then on the unit sphere one draws the great circles perpendicular to $a$ and $b$, oriented according to the direction of the rotation. From an intersection point $O$ of the two great circles one moves backwards along the first great circle by an angle $\frac{\theta}2$, thus reaching a point $A$, and forwards along the second great circle by an angle $\frac{\varphi}2$, thus reaching a point $B$ (cf.~Figure~\ref{Figure2}). Then the composition of the two rotations is a rotation about the axis perpendicular to the great circle through $A$ and $B$, from $A$ to $B$, by an angle $\psi$ equal to twice the angle between $A$ and $B$.
%
%
%

\end{document}